\documentclass{article}
\usepackage{amsmath,amsfonts,amsthm}
\usepackage{graphicx}

\newtheorem{st}{Theorem}
\newtheorem{cor}{Corollary}
\newtheorem{q}{Question}

\title{A cabling formula for the colored Jones polynomial}
\author{Roland van der Veen}
\begin{document}

\maketitle

\begin{abstract}
\noindent We prove an explicit cabling formula for the colored Jones polynomial. As an application we prove the volume conjecture for all zero volume knots and links, i.e. all knots and links that are obtained from the unknot by repeated cabling and connected sum.
\end{abstract}

\section{Introduction}

In this note we study how the unnormalized colored Jones polynomial or quantum $sl_2$ invariant of a link changes under the operation of cabling. We work with a banded link or ribbon link $L$ so that every component is an embedded annulus. Given a diagram $D$ of a banded link inside an annulus we can construct a satellite of $L$ by embedding $D$ into a component $L_i$ of $L$. The $(r,s)$-cabling operation is the special case where we take $D$ to be the closure of the $(r,s)$-torus braid $B_s^r = (\sigma_1\cdots \sigma_{s-1})^r$, where $r\in \mathbb{Z},s\in\mathbb{N}$. To turn $B^r_s$ into a banded tangle we use the blackboard framing and add a positive curl to every overpassing arc, see figure 1 below. The banded link obtained by $(r,s)$-cabling the component $L_i$ of a banded link $L$ will be denoted by $L_{i;s}^r$, we will also call it the $(i;r,s)$-cabling of $L$.

In order to state our cabling formula we need to introduce the following generalizations of the trinomial (not multinomial) coefficients defined in \cite{Andrews}. For a vector $\mathbf{N} = (n_0,\hdots,n_{g-1})$ define ${g \choose w}_{\mathbf{N}}$ to be the coefficient of $x^w$ in the expansion of the product $\prod_{k=0}^{g-1}(x^{\frac{N_k-1}{2}}+x^{\frac{N_k-1}{2}-1}+\hdots + x^{-\frac{N_k-1}{2}})$.

\begin{st} Let $g = \mathrm{gcd}(r,s)$, $p = s/g$ and $\mathbf{N} = (N_0,\cdots,N_{g-1})$. The unnormalized colored Jones polynomial of the $(i;r,s)$-cabling of a banded link $L$ with $c$ components can be expressed as follows:
$J_{M_1,\hdots,M_{i-1},\mathbf{N},\hdots,M_{c}}(L^r_{i;s})(q)=$ \[\sum_{w = -\frac{|\mathbf{N}|-g}{2}}^{\frac{|\mathbf{N}|-g}{2}} {g \choose w}_{\mathbf{N}} q^{\frac{rw(wp+1)}{g}}J_{M_1,\hdots,M_{i-1},2wp+1,\hdots,M_{c}}(L)(q)\]
\end{st}

\noindent In the statement of the theorem we have used the notation $|\mathbf{N}| = N_0+\hdots+N_{g-1}$ and the convention that $J_{M_1,\hdots,M_{i-1},-j,\hdots,M_{c}}(L)(q) = -J_{M_1,\hdots,M_{i-1},j,\hdots,M_{c}}(L)(q)$. 

\begin{figure}[here]
\begin{center}
\includegraphics[]{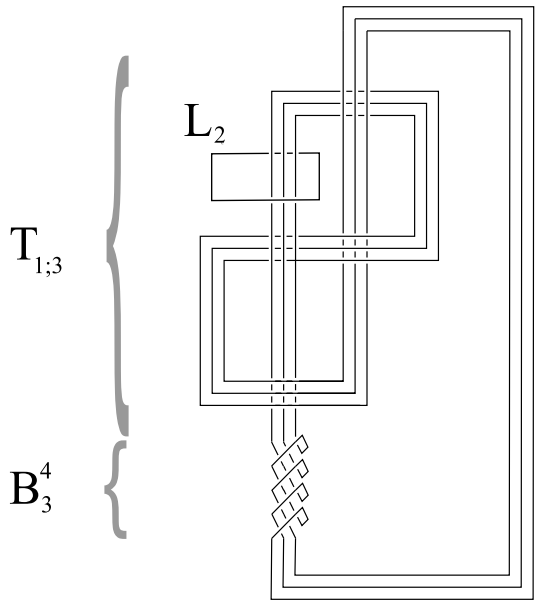}
\caption{We have drawn the link $L_{1;3}^4$, the $(1;3,4)$-cabling of $L = (L_1,L_2)$, where $L_1$ is the figure eight knot and $L_2$ is an unknot. We have indicated the torus braid $B_3^4$ and the opened tangle $T_{1;3}$ mentioned in section 2.} 
\end{center}
\end{figure}

\noindent In the case where $g =1$ and $L$ is the unknot the above cabling formula agrees with Morton's formula for the $(r,s)$-torus knot \cite{Morton}, where his variables are related to ours as $s^2 = q, m = r, p = s$. The case of a $(r,2)$-cabling is also known \cite{ZhengHao}. In all other cases our formula seems to be new.

Our main motivation for proving such a formula is to verify the volume conjecture \cite{Kashaev},\cite{MurakamiMurakami} in the cases where cabling is involved. The volume conjecture states that the normalized colored Jones polynomial $J'(L)$ of a link $L$ determines the simplicial volume of the link complement as follows: \[ \lim_{N\to \infty}\frac{2\pi}{N}|J'_{N,N,\hdots,N}(L)(e^{\frac{2\pi i}{N}})| = \mathrm{Vol}(\mathbb{S}^3-L)\]
\noindent As an immediate corollary to our cabling formula we will find a proof of the volume conjecture for all knots and links whose complement has zero simplicial volume. This is because it is shown in \cite{Gordon} that all such links can be obtained from the unknot by repeated cabling and connected sum. Using the cabling formula we can therefore in principle write down the colored Jones polynomial of any such link. A simple estimate is enough to prove that there cannot be exponential growth, see also section 3. 

\begin{cor}
The volume conjecture is true for all zero volume knots and links.
\end{cor}

\noindent So far the only zero volume and links for which the volume conjecture has been proven are the torus knots \cite{KashTirk} and the $(2,s)$-torus links \cite{ZhengHao}, \cite{Hikami2}.

The cabling formula of Theorem 1 makes it possible to conduct a detailed study of cabled knots and iterated cabling. In the context of the volume conjecture some natural questions would be the following: 

\begin{q}Is the volume conjecture stable under cabling? \end{q}

\noindent Since cabling does not contribute to the simplical volume this would mean that the existing exponential growth is unchanged under cabling.

\begin{q}What is the exact asymptotic expansion for a zero volume link?\end{q}

\noindent In \cite{KashTirk} and \cite{DubKash} explicit asymptotic expansions are given in the case torus knots at the $N$-th root of unity. Away from the root of unity the asymptotics of torus knots has also been studied in the context of the generalized volume conjecture \cite{MurakamiGen}. It would be especially interesting to see whether the polynomial growth related to the roots of the Alexander polynomial predicted and studied in \cite{MurakamiHikami} persists.

\begin{q}How is the colored Jones polynomial of a zero volume link related to q-series identities? \end{q}

\noindent In the case of torus knots and and some special torus links K. Hikami has shown many interesting relations between colored Jones polynomials, modular forms and q-series \cite{Hikami1}, \cite{Hikami2}. 

\begin{q}What is the behavior of the non-commutative A-polynomial of a knot under cabling?\end{q}

\noindent It is known \cite{GaLe} that the colored Jones polynomial satisfies a linear recursion relation. This relation can be encoded in a two variable polynomial with q-coefficients called the non-commutative A-polynomial and it is conjectured \cite{GAJ} to be related to the character variety of the knot group. The knot group behaves well under cabling and according to our cabling formula so does the colored Jones polynomial. It would be interesting to see how these two relate. As before the zero volume knots provide many cases where explicit computations can be made. So far these have only been done for torus knots \cite{Hikami1}.\newline

\noindent An investigation of these questions will be postponed to a subsequent publication. \newline

\noindent \textbf{Acknowledgment} I would like to thank the organizers of the International Conference on Quantum topology 2007 in Hanoi for providing the atmosphere that got this paper started and Stavros Garoufalidis and Hitoshi Murakami for stimulating conversations.

\section{Proof of the cabling formula}

Let us fix a banded link $L$ with $c$ components $L_1,\hdots, L_c$. Choosing a component $L_i$ of $L$ and opening it up we can write $L$ as the closure of a banded $(1,1)$-tangle $T_i$. Define $T_{i;s}$ to be the banded $(s,s)$-tangle obtained from $T_i$ by replacing the opened component $L_i$ of $T_i$ by $s$ parallel bands. In terms of these tangles the $(i;r,s)$-cabling of $L$ is equal to the closure of the composition $T_{i;s} \circ B_s^r$, see also figure 1.

Note that the link $L_{i;s}^r$ has $c+g-1$ components, where $g = \mathrm{gcd}(r,s)$, since the component $L_i$ of $L$ is replaced by the $g$-component torus link that is the closure of $B^r_s$. If we number the strands of the braid $B^r_s$ starting at $0$ then two strands are in the same component of the closure if and only if their numbers are congruent modulo $g$. Therefore each component of the closed braid consists of $p = s/g$ strands. Let us suppose that our link $L_{i;s}^r$ is colored by the integers $(M_1,\hdots,M_{i-1},\mathbf{N},\hdots,M_{c})$, where the vector $\mathbf{N} = (N_0,\hdots,N_{g-1})$ represents the new colors used to color the components of the torus link that replaces the component $L_i$. 

We denote the the $N$-dimensional irreducible representation of quantum $sl_2$ by $V_N$. Using the above coloring, the torus braid corresponds to a morphism from $V = \bigotimes_{k=0}^{s-1} V_{N_{k\ \mathrm{mod}g}}$ to itself. In terms of these morphisms we can now state that the unnormalized $(M_1,\hdots,M_{i-1},\mathbf{N},\hdots,M_c)$-colored Jones polynomial of the $(i;r,s)$-cabling of $L$ (abbreviated by $J$) is the following quantum trace: 

\begin{equation}\label{Jonesquantumtrace}J = J_{M_1,\hdots,M_{i-1},\mathbf{N},\hdots,M_c}(L^r_{i;s}) = \mathrm{Tr}_q (B_s^r \circ T_{i;s},V)\end{equation}

\noindent The first step in calculating this trace is to expand it using the isotypical decomposition of the tensor product as quantum $sl_2$ representations: 

\begin{equation}\label{Isodec}V=\bigotimes_{k=0}^{s-1} V_{N_{k\ \mathrm{mod} g}} \cong \bigoplus_{j=1}^{p|\mathbf{N}|-s+1} \mathrm{Hom}(V_j,V) \otimes V_j\end{equation}

\noindent The range of the summation is calculated by the Clebsch-Gordan rule \cite{Procesi}. The isomorphism maps $(\alpha_j\otimes v)$ on the right to $\alpha_j(v)$ and $B^r_s \circ T_{i;s}$ acts on $\mathrm{Hom}(V_j,V)$ only, while quantum $sl_2$ acts on the $V_j$ only. We can therefore expand the quantum trace in \eqref{Jonesquantumtrace} as follows: 
\begin{equation}\label{Isotrace}J = \mathrm{Tr}_q(B^r_s \circ T_{i;s},V) = \sum_j \mathrm{Tr}(B^r_s \circ T_{i;s},\mathrm{Hom}(V_j,V))\mathrm{Tr}_q(\mathrm{Id},V_j)\end{equation} 
\noindent We know that $\mathrm{Tr}_q(\mathrm{Id},V_j) = [j] = (q^{j/2}-q^{-j/2})/(q^{1/2}-q^{-1/2})$. To calculate the other traces we study how $B^r_s \circ T_{i;s}$ acts on an element $\alpha \in \mathrm{Hom}(V_j,V)$. The action is by composition so first we look at $T_{i;s} \circ \alpha$, the action of $T_{i;s}$. According to the graphical calculus \cite{Turaev} the colored tangle $T_{i;s}$ represents the same map as the $(1,1)$-tangle $T_i$ colored with $V$. Furthermore we can depict $\alpha$ as a coupon connecting the top of a vertical strand colored $V_j$ to the lower end of $T_i$ (this end is colored $V$), see figure 2. 

Now we can slide the coupon $\alpha$ up along the $(1,1)$-tangle $T_{i}$. We obtain $\alpha$ on top of $T_i$, where $T_i$ is now colored with $V_j$. Since $V_j$ is irreducible, the operator represented by $T_i$ is multiplication by a scalar. By closing $T_i$ one sees that this scalar is exactly $[j]^{-1}J_{M_1,\hdots,M_{i-1},j,\hdots,M_c}(L)$. We conclude that $T_{i;s} \circ \alpha = [j]^{-1}J_{M_1,\hdots,M_{i-1},j,\hdots,M_c}(L)\alpha$.

\begin{figure}[here]
\begin{center}
\includegraphics[width = 4 cm,height = 5 cm]{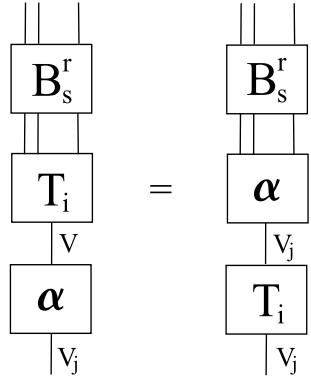}
\caption{Sliding the coupon $\alpha$ up through the tangle $T_{i;s}$.} 
\end{center}
\end{figure}

\noindent The above argument shows that $B^r_s \circ T_{i;s} \circ \alpha = [j]^{-1}J_{M_1,\hdots,M_{i-1},j,\hdots,M_c}(L) B^r_s \circ \alpha$. Therefore equation \eqref{Isotrace} becomes: 
\begin{equation}\label{JIs} J = \sum_{j=1}^{p|\mathbf{N}|-s+1}\mathrm{Tr}(B^r_s,\mathrm{Hom}(V_j,V))J_{M_1,\hdots,M_{i-1},j,\hdots,M_c}(L)\end{equation}

\noindent Note that an equation such as \eqref{JIs} remains valid when one replaces $B^r_s$ by any braid $B$. Such a satellite formula also appears in \cite{MortonStrickland}.

We now concentrate on calculating the trace $\mathrm{Tr}(B_s^r,\mathrm{Hom}(V_j,V))$. Note that $(B_s^r)^s = ((\sigma_1\cdots \sigma_{s-1})^s)^r$ is a central element of the braid group. The element $C = (\sigma_1\cdots \sigma_{s-1})^s$ can also be thought of as a curl in $s$ parallel strands (hence the extra curls in the definition of $B^r_s$), so it acts on $\alpha$ as $C\circ \alpha = q^{\frac{j^2-1}{4}}\alpha$. If we set $D = q^{-\frac{r}{s}\frac{j^2-1}{4}}B^r_s$, then $D^s$ acts as the identity on $\mathrm{Hom}(V_j,V)$. An argument by Jones and Rosso \cite{JonesRosso} shows that its trace does not depend on $q$. For the record we express the trace we were looking for in terms of $D$:
\begin{equation}\label{BD}\mathrm{Tr}(B^r_s,\mathrm{Hom}(V_j,V)) = q^{\frac{r}{s}\frac{j^2-1}{4}}\mathrm{Tr}(D,\mathrm{Hom}(V_j,V))\end{equation}
\noindent Since the trace of $D$ does not depend on $q$ we can calculate it at $q = 1$. In this case $B^r_s = D$ equals the permutation $(\sigma_1\cdots \sigma_{s-1})^r$ and the representation $V$ can be viewed as an $SL(2)$ representation. By collecting common factors in the tensor product we see that $V \cong \bigotimes_{k=0}^{g-1} V_{N_k}^{\otimes p}$, where $p = s/g$. Note that our permutation $D$ acts on $V$ by permuting the factors inside each tensor power. More precisely $D = c_0\cdots c_{g-1}$ where the $c_k$ are disjoint $p$-cycles and $c_k$ permutes the factors of $V_{N_k}^{\otimes p}$. We can therefore interpret $D$ as an element of the Cartesian product $S_p^g$.

One can view each of the tensor powers $V_{N_k}^{\otimes p}$ as a $S_p \times GL(V_{N_k})$ representation, where $S_p$ acts by permuting the tensor factors and $GL(V_{N_k})$ acts diagonally. By Schur-Weyl duality \cite{Procesi} this space allows a simultaneous decomposition into irreducibles:

\[V_{N_k}^{\otimes p} \cong \bigoplus_{\lambda_k \vdash p}E_{\lambda_k}\otimes W_{\lambda_k}\] 
\noindent Here $E_{\lambda_k}$ is an irreducible $S_p$ representation and $W_{\lambda_k}$ is an irreducible $GL(V_{N_k})$ representation and the sum is over all partitions of length no more than $N_k$. Taking the tensor product over all such powers and rearranging the factors gives the following decomposition of $V$ as an $S_p^g\times GL(V_{N_0})\times \cdots \times GL(V_{N_{g-1}})$ representation.

\[V \cong \bigoplus_{\lambda_{0}\hdots \lambda_{g-1} \vdash p}E_{\lambda_{0}}\otimes \cdots \otimes E_{\lambda_{g-1}}\otimes  W_{\lambda_0}\otimes \cdots \otimes W_{\lambda_{g-1}}\]

\noindent The $GL(V_{N_0})\times \cdots \times GL(V_{N_{g-1}})$ representation $W_{\lambda_0}\otimes \cdots \otimes W_{\lambda_{g-1}}$ is in a natural way also an $SL(2)$-representation. We can therefore decompose it into irreducible $SL(2)$-representation as follows:

\[W_{\lambda_0}\otimes \cdots \otimes W_{\lambda_{g-1}} \cong \bigoplus_{j}\mathrm{Hom}(V_j,W_{\lambda_0}\otimes \cdots \otimes W_{\lambda_{g-1}})\otimes V_j\] 

\noindent Hence we find the following decomposition of $V$:

\[V \cong \bigoplus_{j}\left(\bigoplus_{\lambda_0\hdots \lambda_{g-1} \vdash p}E_{\lambda_0}\otimes \cdots \otimes E_{\lambda_{g-1}}\otimes  \mathrm{Hom}(V_j,W_{\lambda_0}\otimes \cdots \otimes W_{\lambda_{g-1}})\right)\otimes V_j\]

\noindent Note that we only allow $\lambda_k$ to have length $\leq N_k$. Comparing this decomposition with the isotypical decomposition \eqref{Isodec} for $q = 1$ we see that:

\begin{equation}\label{Homdec}\mathrm{Hom}(V_j,V) \cong \bigoplus_{\lambda_0\hdots \lambda_{g-1} \vdash p}E_{\lambda_0}\otimes \cdots \otimes E_{\lambda_{g-1}}\otimes  \mathrm{Hom}(V_j,W_{\lambda_0}\otimes \cdots \otimes W_{\lambda_{g-1}})\end{equation}

\noindent The action of $D$ on this space is only on the $E_{\lambda_k}$. As above we write $D  = c_0\cdots c_{g-1}$ as a product of disjoint $p$-cycles and we note that the cycle $c_k$ acts in $E_{\lambda_k}$. If we define $R_{j,\lambda_0\hdots,\lambda_{g-1}}=\mathrm{dim}(\mathrm{Hom}(V_j,W_{\lambda_0}\otimes \cdots \otimes W_{\lambda_{g-1}}))$ and denote all $c_k$ by $c$ we have the following expression for the trace: 
\begin{equation}\label{TrD}\mathrm{Tr}(D,\mathrm{Hom}(V_j,V)) =  \sum_{\lambda_0,\hdots \lambda_{g-1} \vdash p} \chi_{\lambda_0}(c)\cdots\chi_{\lambda_{g-1}}(c)\ R_{j,\lambda_0\hdots,\lambda_{g-1}}\end{equation}

\noindent Here $\chi_{\lambda}$ is the character of the symmetric group $S_p$. To calculate these traces we note that the dimensions $R_{j,\lambda_1\hdots,\lambda_g}$ are encoded in the following product expansion of Schur functions \cite{Macdonald}:
\begin{equation}\label{GenSchur}\prod_{k=0}^{g-1} s_{\lambda_k}(t^{N_k-1},t^{N_k-3},...,t^{-N_k+1}) = \sum_{\mu \vdash p|\mathbf{N}|-s}R_{\mu_1-\mu_2+1,\lambda_0\hdots,\lambda_{g-1}}s_{\mu}(t,t^{-1})\end{equation}

\noindent The above condition on the direct sum in \eqref{Homdec} that the length of $\lambda_k$ is no more than $N_k$ is now incorporated into the Schur functions, because an $N$-variable Schur function $s_{\lambda}$ is nonzero only if $\lambda$ has no more than $N$ parts.

As a quick check that the above formula \eqref{GenSchur} includes all the dimensions $R$, note that $\lambda_k \vdash p$ and so the highest degree term of $s_{\lambda}(t^{N_k-1},t^{N_k-3},...,t^{-N_k+1})$ is $t^{p N_k-p}$, hence $\mu \vdash p|\mathbf{N}|-s$ because the highest degree term of $s_{\mu}(t,t^{-1})$ is $t^{\mu_1}$. Note also that only length-two partitions $\mu$ contribute to the sum because the Schur functions $s_{\mu}$ have only two variables. In case $\mu = (\mu_1 \mu_2)$ we have $s_{\mu}(x,y) = \frac{x^{\mu_1+1}y^{\mu_2}-x^{\mu_2}y^{\mu_1+1}}{x-y}$ and $s_{\mu}(t,t^{-1}) = \frac{t^{\mu_1+1-\mu_2}-t^{-(\mu_1+1-\mu_2)}}{t-t^{-1}}$. The range of the summation in \eqref{GenSchur} therefore agrees exactly with the range of $j$ in the decomposition of the tensor product \eqref{Isodec}. From now on we will parametrize this range by two part partitions $\mu \vdash p|\mathbf{N}|-s$ and set $V_\mu = V_j$, where $j = \mu_1-\mu_2+1$.

Using the product of Schur functions in \eqref{GenSchur} as a generating function for the dimensions $R$ we can compute $\mathrm{Tr}(D,\mathrm{Hom}(V_{\mu},V))$ simultaneously for all $\mu$. First we sum the right hand side and the left hand side of equation \eqref{GenSchur} with the characters as in the sum in \eqref{TrD}: 

\begin{equation}\label{GenR}\sum_{\mu \vdash |\mathbf{N}|-s}\left(\sum_{\lambda_0,\hdots \lambda_{g-1} \vdash p} \chi_{\lambda_0}(c)\cdots\chi_{\lambda_{g-1}}(c)R_{\mu,\lambda_0\hdots,\lambda_{g-1}}\right)s_{\mu}(t,t^{-1}) =\end{equation}
\[\sum_{\lambda_0,\hdots\lambda_{g-1} \vdash p} \chi_{\lambda_0}(c)\cdots \chi_{\lambda_{g-1}}(c)\prod_{k=0}^{g-1} s_{\lambda_k}(t^{N_k-1},t^{N_k-3},...,t^{-N_k+1}) = \]
\[\prod_{k=0}^{g-1} \sum_{\lambda_k\vdash p}\chi_{\lambda_k}(c)s_{\lambda_k}(t^{N_k-1},t^{N_k-3},...,t^{-N_k+1}) = \]
\begin{equation}\label{PowProd}\prod_{k=0}^{g-1}p_p(t^{N_k-1},t^{N_k-3},...,t^{-N_k+1})\end{equation}

\noindent In the last line we used the expansion of the power sum function $p_p$ in terms of Schur functions \cite{Macdonald} p.114. We can use the notion of plethysm \cite{Macdonald} p.138 to write the power sums in terms of the complete symmetric functions $h$ as follows: \[p_p(t^{N_k-1},t^{N_k-3},...,t^{-N_k+1}) = (p_p \circ h_{N_k-1})(t,t^{-1}) = h_{N_k-1}(t^p,t^{-p})\]

\noindent As with partitions we define $h_{\mathbf{N}-1} = \prod_{k=0}^{g-1}h_{N_k-1}$ so that we can write the above product \eqref{PowProd} as:

\begin{equation}\label{HAnswer}\prod_{k=0}^{g-1}p_p(t^{N_k-1},t^{N_k-3},...,t^{-N_k+1})=h_{\mathbf{N}-1}(t^p,t^{-p})\end{equation}

\noindent Comparing equations \eqref{GenR} and \eqref{HAnswer} we see that to compute the traces in \eqref{TrD} we need to express $h_{\mathbf{N}-1}(t^p,t^{-p})$ in terms of the $s_{\mu}(t,t^{-1})$. The coefficients in this expansion will be our traces. 

First we consider the case where $p=1$. Since we have $s = g$ there are $\left\lfloor (|\mathbf{N}|-g)/2\right\rfloor+1$ two part partitions $\mu \vdash |\mathbf{N}|-s$, that we will number by the half integers $w = (\mu_1-\mu_2)/2$. Since $h_{\mathbf{N}-1}(t,t^{-1})$ is a product of geometric series we can use the notation defined in the introduction to write: $h_{\mathbf{N}-1}(t,t^{-1}) = $
\begin{equation}\label{P1geval}\sum_{w=(|\mathbf{N}|-g)/2-\left\lfloor (|\mathbf{N}|-g)/2\right\rfloor}^{(|\mathbf{N}|-g)/2}\left({g \choose w}_{\mathbf{N}}- {g \choose w+1}_{\mathbf{N}}\right)\frac{t^{2w+1}-t^{-2w-1}}{t-t^{-1}}\end{equation}

\noindent The above equation settles the case $p=1$ since the Schur functions are exactly the geometric series in the two variable case. We will come back to the $p=1$ case after we have found the coefficients for $p>1$. We will then show that these cases can be unified into a single cabling formula for the colored Jones polynomial.

From now on let's assume that $p\geq 2$. We will carry out the expansion of $h_{\mathbf{N}-1}(t^p,t^{-p})$ into two variable Schur functions in a little more generality by setting $x_1 = t$ and $x_2 = t^{-1}$. As a first step consider the following fundamental identity \cite{Macdonald} p.62 relating the complete symmetric functions to the Schur functions:

\begin{equation}\label{FundId} \sum_{\lambda}h_{\lambda}(x)m_{\lambda}(y) = \prod_{i,j}(1-x_iy_j)^{-1}  = \sum_{\lambda}s_{\lambda}(x)s_{\lambda}(y)\end{equation}

\noindent Here $m_{\lambda}(y)$ is the monomial symmetric function and $x = (x_1,x_2,\hdots)$ and the sums range over all partitions $\lambda$. The first part of the identity \eqref{FundId} implies that $h_{\mathbf{N}-1}(x_1^p,x_2^p)$ is the coefficient of $\mathbf{y}^{p(\mathbf{N}-1)} = y_0^{p (N_0-1)}\cdots y_{g-1}^{p (N_{g-1}-1)}$ in the expansion of the product $\prod_{k=0}^{g-1}(1-x_1^py_k^p)^{-1}(1-x_2^py_k^p)^{-1}$. Factoring we can also write this product as $\prod_{k=0}^{g-1}\prod_{u=0}^{p-1}(1-x_1\omega^u y_k)^{-1}(1-x_2\omega^u y_k)^{-1}$, where $\omega = \exp(2\pi i/p)$. Using the second half of the identity \eqref{FundId} we can expand the same product into two variable Schur functions: 
\[\prod_{k=0}^{g-1}\prod_{u=0}^{p-1}(1-x_1\omega^u y_k)^{-1}(1-x_2\omega^u y_k)^{-1} = \sum_{\mu}s_{\mu}(y_0,\omega y_0,\hdots,\omega^{p-1} y_{g-1})s_{\mu}(x_1,x_2)\]

\noindent Taking the coefficient of $\mathbf{y}^{p(\mathbf{N}-1)}$ on the left gives $h_{\mathbf{N}-1}(x_1^p,x_2^p)$ back, while the coefficient of $\mathbf{y}^{p(\mathbf{N}-1)}$ on the right hand side is $\sum_\mu K_{\mu}s_{\mu}(x_1,x_2)$, where $K_\mu$ is the coefficient of $\mathbf{y}^{p(\mathbf{N}-1)}$ in $s_{\mu}(y_0,\omega y_0,\hdots,\omega^{p-1} y_{g-1})$. It follows that the trace we wanted to calculate is equal to this coefficient: $\mathrm{Tr}(D,\mathrm{Hom}(V_{\mu},V)) =  K_{\mu}$.

We proceed with the calculation of $K_{\mu}$. By \cite{Macdonald} p.72 we can expand the Schur function $s_{\mu}(y_0,\omega y_0,\hdots,\omega^{p-1} y_{g-1})$ in terms of skew Schur functions making it possible to extract the powers of $y_j$.
\[s_{\mu}(y_0,\omega y_0,\hdots,\omega^{p-1} y_{g-1}) = \sum \prod_{j = 0}^{g-1} y_j^{|\nu^{(j)}-\nu^{(j-1)}|}s_{\nu^{(j)}/\nu^{(j-1)}}(1,\omega,\hdots,\omega^{p-1})\]

\noindent The sum ranges over all sequences of partitions $0 = \nu^{(0)}\subset \nu^{(1)} \hdots \subset \nu^{(g)} = \mu$. Only the sequences of partitions satisfying $|\nu^{(j)}-\nu^{(j-1)}| = p(N_j-1)$ contribute to the coefficient of $\mathbf{y}^{p(\mathbf{N}-1)}$ in this sum. The contribution of such a sequence is $\prod_{j = 0}^{g-1} s_{\nu^{(j)}/\nu^{(j-1)}}(1,\omega,\hdots,\omega^{p-1})$. In the cases we are interested in $\mu$ has length $2$ and it follows from \cite{Macdonald} p.91 that this contribution is zero unless $\nu^{(j)}$ can be obtained from $\nu^{(j-1)}$ by attaching $N_j-1$ border strips of length $p$. A set theoretical difference between two nested partitions is called a border strip of length $p$ if it is connected, does not contain any $2\times 2$ squares and has $p$ elements. If this is the case for all $j$, then the contribution is $\sigma_p(\mu) = (-1)^{r_2}$, where $0\leq r_2<p$ is the residue of $\mu_2$ modulo $p$.

The above argument shows that $K_{\mu} = \sigma_p(\mu)|K_{\mu}|$, and that $|K_{\mu}|$ is the number of sequences of partitions $0 = \nu^{(0)}\subset \nu^{(1)} \hdots \subset \nu^{(g)} = \mu$ satisfying the condition that $\nu^{(j)}$ is obtained from $\nu^{(j-1)}$ by attaching $N_j-1$ border strips of length $p$. 

To count the number of such sequences we first look a all ways of constructing two part partitions by attaching length $p$-border strips, starting with the empty partition. One notices that the lower leftmost corner of a new border strip can be attached at exactly two places: either the end of the first row or the end of the second row. It follows that we can list all possibilities in a Pascal-like triangle, see figure 3. 

\begin{figure}[here]
\begin{center}
\includegraphics[width = 12 cm,height = 5 cm]{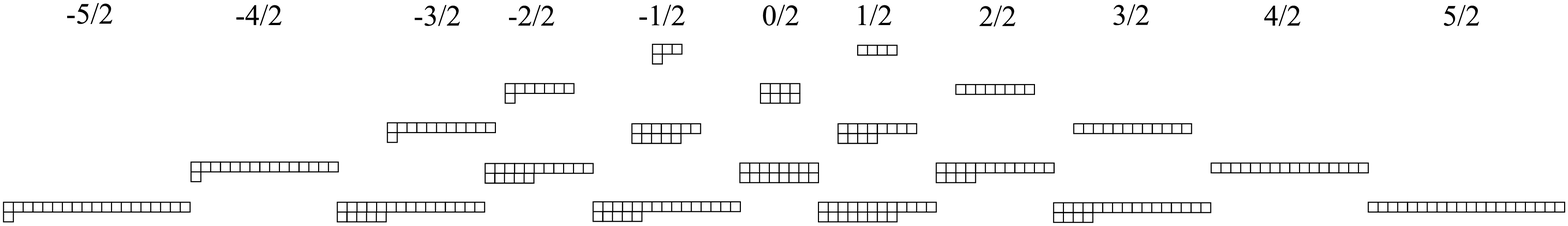}
\caption{The Pascal-like triangle of all two part partitions that can be glued from $5$ or less border strips of length $p = 4$.} 
\end{center}
\end{figure}

\noindent The $k$-th row of our Pascal triangle contains $k+1$ partitions that we index by the half integers $\frac{-k}{2},\frac{-k}{2}+1,\hdots,\frac{k}{2}$. The triangle is ordered such that if $\mu$ has index $w(\mu) = w$ then $\mu_1-\mu_2+1 = \mathrm{sgn}(w)(2wp+1)$ and $\mathrm{sgn}(w) = \sigma_p(\mu)$, provided that we agree that $\mathrm{sgn}(0)=1$.

Interpreting $|K_{\mu}|$ in terms of the Pascal triangle we see that $|K_{\mu}|$ is equal to the number of ways to move from the (empty) top of the triangle to $\mu \vdash p|\mathbf{N}|-s$ in $g$ steps of length ${N_k}-1$, for $k = 0\hdots,g-1$. A step of length $A$ can go in $A+1$ directions that are conveniently indexed by $\frac{-A}{2},\frac{-A}{2}+1,\hdots,\frac{A}{2}$. The sum of these indices must be the index of $\mu$, that is $w$. Therefore the number of ways is exactly the coefficient of $x^w$ in the product $\prod_{k=0}^{g-1}(x^{-\frac{N_k-1}{2}}+x^{-\frac{N_k-1}{2}+1}+\hdots+x^{\frac{N_k-1}{2}})$. Using the notation defined in the introduction we see that $|K_{\mu}| = {g \choose w}_{\mathbf{N}}$ and so $K_{\mu} = \mathrm{sgn}(w){g \choose w}_{\mathbf{N}}$. 

We can now summarize the calculation of $J = J_{M_1,\hdots,M_{i-1},\mathbf{N},\hdots,M_c}(L^r_{i;s})$ as follows. First we expressed it in equation \eqref{JIs} as a sum over the colors $j$ of a trace times the colored Jones of the uncabled link. In \eqref{BD} the dependence on $q$ was extracted and finally the rest of the trace was shown to be given by the coefficients in equation \eqref{P1geval} for $p=1$ and by $K_\mu$ for $p\geq 2$. Along the we changed the sum over colors into a sum over two part partitions $\mu$, via $j = \mu_1-\mu_2+1$. We just saw that we can enumerate the $\mu$ for which $K_\mu$ is nonzero by the half integers $w$. Combining all this we find for $p\geq 2$: $J_{M_1,\hdots,M_{i-1},\mathbf{N},\hdots,M_c}(L^r_{i;s})=$
\[\sum_{w=-(|\mathbf{N}|-g)/2}^{(|\mathbf{N}|-g)/2}q^{\frac{rwp(wp+1)}{s}}\mathrm{sgn}(w){g \choose w}_{\mathbf{N}}J_{M_1,\hdots,M_{i-1},|2wp+1|,\hdots,M_c}(L)\]

\noindent Remarkably this formula is also valid for $p=1$. In this case the non-negative $w$ run through all partitions and both $w$ and $-(w+1)$ contribute to the same term $J_{M_1,\hdots,M_{i-1},|2w+1|,\hdots,M_c}(L)$, thus producing the coefficients found in \eqref{P1geval}.
 
Finally the convention $J_{M_1,\hdots,M_{i-1},-k,\hdots,M_c} = -J_{M_1,\hdots,M_{i-1},k,\hdots,M_c}$ allows us to incorporate the sign into the colored Jones polynomial because $\mathrm{sgn}(w) = \mathrm{sgn}(2wp+1)$. Noting that $p = s/g$ the proof of theorem 1 is complete.

\section{Proof of the volume conjecture for zero volume links}

As an application of the cabling formula we just proved we will show that the volume conjecture holds true for all zero volume knots and links (corollary 1). Such links are related to cabling by the fact that all links of zero simplicial volume are obtained from the unknot by repeated cabling and connected sum \cite{Gordon}. 

In the volume conjecture the usual normalization of the colored Jones polynomial for knots is to divide by the value of the unknot. However, if a link has $s$ split components then the unnormalized Jones has at least an $s$ fold zero at the $N$-th root of unity \cite{vanderVeen}. We therefore choose to normalize it by dividing by $[N]^s$. For convenience we work with $A = q^{1/4}$ instead of $q$ so that for the volume conjecture we need to evaluate at $A = e^{\pi i/2N}$.

The unnormalized colored Jones polynomial of a link with $c$ components is a multi-sequence of Laurent polynomials in $A$ of the form $P_{\mathbf{N}}(A) = \sum_j a_{j,\mathbf{N}}A^j$, where $\mathbf{N} = (N_1,\hdots,N_c)$. We say that a multi-sequence is moderate if there exist $C,m>0$ such that  \[ \forall\ \mathbf{N}:\quad |\mathrm{maxdeg}P_{\mathbf{N}}|,|\mathrm{mindeg}P_{\mathbf{N}}|,|a_{j,\mathbf{N}}| < C|\mathbf{N}|^m\] \noindent It is clear that the set of moderate sequences is closed under products and differentiation with respect to $A$.

To prove the volume conjecture for all zero volume links we first show that their unnormalized colored Jones polynomials are moderate. Since such links are obtained from the unknot by repeated cabling and connected sum, we need to show that $[N]$ is moderate (trivial) and that the set of moderate multi-sequences is closed under taking products and under an operation that generalizes cabling of links that we now define. Given two multi-sequences $P_{\mathbf{M}}(A)$ and $Q_{\mathbf{N},w}(A)$ define a new multi-sequence $R$ by \[R_{M_1,\hdots M_{i-1}, \mathbf{N},\hdots,M_c}(A) = \sum_{|w|\leq (|\mathbf{N}|-g)/2}Q_{\mathbf{N},2w}(A)P_{M_1,\hdots,M_{i-1},|2wp+1|,\hdots, M_c}(A)\]

\noindent The triangle inequality shows that that if $P$ and $Q$ are moderate then so is $R$. The cabling formula applied to the colored Jones polynomial $J$ of any link is an example of the above defined operation where $P = J$ and $Q_{\mathbf{N},w}(A)$ is a multi-sequence that does not depend on the link and is readily read off from the cabling formula or rather the version at the end of the last section. It is not hard to see that this $Q$ is indeed moderate. 

Now that we know that the unnormalized colored Jones polynomial of a zero volume knot $J_{\mathbf{N}}(A)$ is moderate we conclude the proof of corollary 1 as follows. To evaluate the normalized Jones polynomial we will use 'l Hospital's rule after writing it as \[J_{N,\hdots,N}(A)/[N]^s = \tilde{J}_{N,\hdots,N}(A)/(A^{2N}-A^{-2N})^s\] \noindent where $\tilde{J}$ = $(A^{2}-A^{-2})^sJ$ is again moderate. If we differentiate both numerator and denominator $s$ times with respect to $A$ and evaluate at $e^{\pi i/2N}$ we get $(\frac{e^{\pi i/2N}}{-4N})\frac{d^s}{dA^s}\tilde{J}_{N,\hdots,N}(e^{\pi i/2N})$. Since $\frac{d^s}{dA^s}\tilde{J}_{\mathbf{N}}$ is again moderate, it follows that:
\[\lim_{N\to \infty}\frac{2\pi}{N}|\frac{J_{N,\hdots,N}}{[N]^k}(e^{\frac{\pi i}{2N}})| = 0\]

\end{document}